\documentclass[12pt,leqno]{amsart}
\usepackage{amsmath}
\usepackage{amssymb}

\def\underset#1#2{{\mathrel{\mathop {{}_{} {#2}}\limits_{{#1}_{}}}}}
\def\upplim_#1{\underset{#1}{\overline\lim}\;}
\def\lowlim_#1{\underset{#1}{\underline\lim}\;}
\newtheorem{definition}[equation]{\indent{\it Definition}\rm }

\newtheorem{lem}[equation]{Lemma}
\newtheorem{lemma}[equation]{Lemma}

\newtheorem{theorem}[equation]{Theorem}

\newcommand{\C}{{\mathbf{C}}}
\newcommand{\pnc}{{\mathbf{P}^n{(\mathbf{C})}}}

\renewcommand{\P}{{\mathbf{P}}}

\newcommand{\zero}{\mathrm{Zero}}

\newcommand{\supp}{\mathrm{Supp}\,}

\newcommand{\Z}{\mathbf{Z}}

\numberwithin{equation}{section}

\title[Multiple values and finiteness problem of meromorphic mappings]{Multiple values and finiteness problem of meromorphic mappings sharing different families of moving hyperplanes}

\date { }
\author{Ha Huong Giang}

\address{Faculty of Fundamental Sciences, Electric Power University\\
235-Hoang Quoc Viet, Tu Liem , Ha Noi, Vietnam.}

\email{hhgiang79@yahoo.com}

\begin{document}

\begin{abstract}
In this article, we show some uniqueness theorems for meromorphic mappings of $\C^n$ into the complex projective space  $\pnc$ sharing different families of moving hyperplanes regardless of multiplicites, where all intersecting points between these mappings and moving hyperplanes with multiplicities more than a certain number do not need to be counted.
\end{abstract}

\def\thefootnote{\empty}
\footnotetext{}

\maketitle

\section{Introduction}
In 1926, Nevanlinna \cite{N} showed that for two nonconstant meromorphic functions $f$ and $g$ on the complex plane $\C$, if they have the same inverse images for five distinct values, then $f \equiv g$. After that, many mathematicians have generalized the Nevanlinna's result to the case of meromorphic mappings of $\C^m$ into $\pnc$. Specially, in 1975, Fujimoto \cite{Fu} proved that for two linearly nondegenerate meromorphic mappings $f$ and $g$ of $\C^m$ into $\pnc$, if they have the same inverse images counted with multiplicities for $3n+2$ hyperplanes in general position in $\pnc$, then $f \equiv g$.

In 1983, L.Smiley \cite{Sm} considered meromorphic mappings with share $3n+2$ hyperplanes of $\pnc$ without counting multiplicities and he proved the following.
\vskip0.2cm
\noindent
{\bf Theorem A} (see  \cite{Sm}).\ {\it Let $f,g: \C^m \rightarrow \pnc$ be linearly  nondegenerate meromorphic mappings of $\C^m$ into $\pnc$ . Let $\{H_i\}_{i=1}^q$ $(q \ge 3n+2)$ be hyperlanes in $\pnc$ in general position. Assume that
\begin{align*}
&(i) f^{-1}(H_i) =g^{-1}(H_i),   \quad for \quad  1 \leq i \leq q\\
&(ii) \dim(f^{-1}(H_i)\cap f^{-1}(H_j)) \leq m - 2, \quad   for \quad all \quad  1 \leq i < j \leq q\\
&(iii) f = g \quad on \quad  \bigcup_{i=1}^qf^{-1}(H_i)
\end{align*}
then $f = g$.}

\vskip0.2cm
In 2010, Gerd Dethloff , Si Duc Quang and Tran Van Tan \cite{D-Q-T} considered the case where the mappings sharing different families of hyperplanes. They showed that

\vskip0.2cm
\noindent
{\bf Theorem B} (see  \cite{D-Q-T}).\ {\it Let $f,g: \C^m \rightarrow \pnc$ be a meromorphic mapping. Let $\{H_i\}_{i=1}^q$ and $\{L_i\}_{i=1}^q$, $ (q \geq 3n + 2)$ be families of hyperplanes in $\pnc$ in general position. Assume that
\begin{align*}
&(i) f^{-1}(H_i) =g^{-1}(L_i),   \quad for \quad  1 \leq i \leq q\\
&(ii) \dim(f^{-1}(H_i)\cap f^{-1}(H_j)) \leq m - 2, \quad   for \quad all \quad  1 \leq i < j \leq q\\
&(iii) \frac{(f,H_i)}{(g,L_i)} = \frac{(f,H_j)}{(g,L_j)}\quad on \quad  \bigcup_{k=1}^qf^{-1}(H_k)\setminus (f^{-1}(H_i)\cap f^{-1}(H_j)),\quad for \quad all \quad  1 \leq i < j \leq q.
\end{align*}
Then the following assertions hold:
$$
dim\langle Imf\rangle=dim\langle Img\rangle:=p
$$
where for a subset $X\subset \pnc$, we denote by $\langle X\rangle$ the smallest projective subspace of $\pnc$ containing $X$.
$$
If \quad q>\dfrac{2n+3-p+\sqrt{(2n+3-p)^2+8(p-1)(2n-p+1)}}{2}(\geq 2n+2)
$$
then
\begin{align*}
\frac{(f,H_1)}{(g,L_1)} \equiv ... \equiv \frac{(f,H_q)}{(g,L_q)}
\end{align*}
Furthermore, there exists a linear projective transformation $\mathcal{L}$ of $\pnc$ into itself such that $\mathcal{L}(f) \equiv g$ and $\mathcal{L}(H_i \cap \langle Imf\rangle)=\mathcal{L}_i \cap \mathcal{L}(\langle Imf\rangle)$ for all $i \in \{1,...,q\}$.}

\vskip0.2cm
In 2011, Ting-Bin Cao and Hong-Xun Yi \cite{C-Y} showed the following result

\vskip0.2cm
\noindent
{\bf Theorem C} (see  \cite{C-Y}).\ {\it Let $f$ and $g$ be two linearly non-degenerate meromorphic mappings of $\C^m$ into $\pnc$, and let $H_1, H_2,...,H_q$ be $q$ $(q\geq 2n)$ hyperplanes in general position such that $dimf^{-1}(H_i\cap H_j)\leq m-2$ for $i \neq j$. Take $m_j$ $(j=1,2,...,q)$ be positive integers or $\infty$ such that $m_1\geq m_2\geq...\geq m_n\geq n$,
 $$
 \nu_{(f,H_j),\leq m_j}^1= \nu_{(g,H_j),\leq m_j}^1 \quad (j=1,2,...,q)
 $$
 and $f(z)=g(z)$ on $\bigcup_{j=1}^q\{z \in \C^m:0<\nu_{(f,H_j)}\leq m_j\}$. If
 $$
 \sum_{j=3}^q\dfrac{m_j}{1+m_j}>\dfrac{nq-q+n+1}{n}-\dfrac{4n-4}{q+2n-2}+\biggl(\dfrac{1}{1+m_1}+\dfrac{1}{1+m_2}\biggl)
 $$
  then $f(z) \equiv g(z)$.}

  Recently, Zhonghua Wang and Zhenhan Tu proved a uniqueness theorem for meromorphic mappings in several complex variables into the complex projective space $\pnc$ with two families of moving targets as follows.
\vskip0.2cm
\noindent
{\bf Theorem D} (see  \cite{W-T}).\ {\it Let $f, g, a_i, b_i: \mathbf{C}^m \rightarrow \pnc$ be meromorphic mappings $(i = 1, 2,...,q)$. Suppose that $\{a_i\}_{i=1}^q$ are ``small'' (with respect to f) and located in the general position, and that $\{b_i\}_{i=1}^q$ are ``small'' (with respect to g) and located in the general position such that f and g are linearly nondegenerate over $\mathcal{R}(\{a_i,b_i\}_{i=1}^q)$. For any reduced representations $a_i = (a_{i0}, ...,a_{in})$ and $b_i = (b_{i0}, ...,b_{in})$ $ (i = 1, 2,...,q)$, we may assume $a_{i0}\not\equiv 0$ and  $b_{i0}\not\equiv 0 $ $(i = 1, 2,...,q)$ by changing the homogeneous coodinate system of $\pnc$. Let $\widetilde{a}_i = \frac{a_i}{a_{i0}} $ and $\widetilde{b}_i = \frac{b_i}{b_{i0}}$ $  (i = 1, 2,...,q)$. Assume that
\begin{align*}
&(i) \nu_{(f,\widetilde{a}_i)}^1(z) = \nu_{(g,\widetilde{b}_i)}^1(z),   for  1 \leq i \leq q\\
&(ii) \dim \{z \in \C^m: (f(z), a_i(z)) = (f(z), a_j(z)) = 0\} \leq m - 2,   for  1 \leq i < j \leq q\\
&(iii) \frac{(f,\widetilde{a}_i)}{(g,\widetilde{b}_i)} = \frac{(f,\widetilde{a}_j)}{(g,\widetilde{b}_j)}  on  \bigcup_{k=1;k\neq i,j}^q\{z \in C^n: (f(z),a_k(z)) = 0\},  for  1 \leq i < j \leq q.
\end{align*}
 Then \\
If $q = 2n^2 + 2n + 3$ then there exist $\{i_1, ...,i_{n+1}\} \subset \{1, ...,q\}$such that
\begin{align*}
\frac{(f,\widetilde{a}_{i_1})}{(g,\widetilde{b}_{i_1})} \equiv ... \equiv \frac{(f,\widetilde{a}_{i_{n+1}})}{(g,\widetilde{b}_{i_{n+1}})}
\end{align*}
which immediately means that there exists a matrix L with its elements $L_{ij}$ in $\mathcal{R}(\{a_i,b_i\})_{i=1}^q$ such that $L(f) = g$}

\vskip0.2cm
We note that, in the above theorem the mappings are assumed to be linearly nondegenerate. Our purpose in this paper is to study the case where the mappings may be degenerate. We will show some uniqueness theorems for mappings sharing different families of moving hyperplanes regardless of multiplicities, which are improvements and extensions of some recent results in this direction when reduced to the case of mappings sharing the same family of moving hyperplanes. Our main results of this work are stated as follows.

\noindent
Let $f^t:\C^m \to \P^n(\C)$ be  meromorphic mapping. Let $\{a_i^t\}_{i=1}^{q}$ be family of  moving hyperplanes in $\P^n(\C)$ in general position such that $a_i^t$ be ``slowly'' with respect to $f^t$. By changing the homogeneous coodinate system of $\pnc$ if necessary, we may assume that  $a^t_{i0} \not \equiv 0$ $(1\le i\le q)$ any given meromorphic mapping $a_i^t=(a^t_{i0},...,a^t_{in})$. Let $\widetilde{a}^t_i=\frac{a^t_i}{a^t_{i0}}$, $1\le i \le q$.

We will prove the following.

\begin{theorem}\label{1.1}
Let $f^1,f^2:\C^m \to \P^n(\C)$ be two meromorphic mappings.  Let $k_i$ $(1\le i\le q)$ be positive integers or $\infty$. Let $\{a_i^t\}_{i=1}^{q}$ $(t=1,2)$ be two families of  moving hyperplanes in $\P^n(\C)$ in general position such that $a_i^t$ is "slowly'' with respect to $f^t$ and $\dim\ \{z\in\C^m: \nu_{(f^t,a^t_i),\le k_i}.\nu_{(f^t,a^t_j),\le k_j}>0\} \leq m-2$ $ (1\le i<j\le q, t=1,2)$. We assum that:

(a) $\min\{\nu_{(f^2,\widetilde{a}_i^2),\leq k_i}(z),1\}=\min\{\nu_{(f^1,\widetilde{a}_i^1),\leq k_i}(z),1\}$ $ (1 \le i \le q),$ for all $z\in\C^m$,

(b)  $\frac{(f^1,\widetilde{a}_i^1)}{(f^2,\widetilde{a}_i^2)} = \frac{(f^1,\widetilde{a}_j^1)}{(f^2,\widetilde{a}_j^2)}$ on $ \bigcup_{\underset{ v\neq i,j}{v=1}}^q\supp\{z \in \C^m:\nu_{(f^1,a^1_v),\le k_v}(z) \}$, for $  1 \leq i < j \leq q.$

If $q > 3n^2+n+2$ and $\sum_{i=1}^q\frac{1}{k_i+1}< \biggl(\frac{2q}{3n(n + 1)}-\frac{2q}{q+2n-2}\biggl)$, then there exist $n+1$ indices $1\le i_1<i_2<\cdots <i_{n+1}\le q$ such that
\begin{align}\label{new2}
\dfrac{(f^1,\widetilde{a}_{i_1}^1)}{(f^2,\widetilde{a}_{i_1}^2)} =\cdots = \dfrac{(f^1,\widetilde{a}_{i_{n+1}}^1)}{(f^2,\widetilde{a}_{i_{n+1}}^2)}.
\end{align}
\end{theorem}

\begin{theorem}\label{1.2}
Let $f^1,f^2,f^3:\C^m \to \P^n(\C)$ be three meromorphic mappings.  Let $k_i$ $(1\le i\le q)$ be positive integers or $\infty$. Let $\{a_i^t\}_{i=1}^{q}$ $(t=1,2,3)$ be three families of  moving hyperplanes in $\P^n(\C)$ in general position such that $a_i^t$ be ``slowly'' with respect to $f^t$ and $\dim\ \{z\in\C^m: \nu_{(f^t,a^t_i),\le k_i}.\nu_{(f^t,a^t_j),\le k_j}>0\} \leq m-2\quad (1\le i<j\le q, 1\le t \le 3)$. We assume that:
\begin{enumerate}
\item[(a)] $\min\{\nu_{(f^t,\widetilde{a}_i^t),\leq k_i}(z),1\}=\min\{\nu_{(f^1,\widetilde{a}_i^1),\leq k_i}(z),1\}$ $(1 \le i \le q, t=2,3),\ \forall z\in\C^m$,
\item[(b)] $\frac{(f^1,\widetilde{a}_i^1)}{(f^t,\widetilde{a}_i^t)} = \frac{(f^1,\widetilde{a}_j^1)}{(f^t,\widetilde{a}_j^t)}$ on $\bigcup_{\underset{ v\neq i,j}{v=1}}^q\supp\{z \in \C^m:\nu_{(f^1,a^1_v),\le k_v}(z) \}$, $1 \leq i < j \leq q,$ $ t=2,3.$
\end{enumerate}
If $q> \dfrac{9n^2+7n+6}{4}$ and $\sum_{i=1}^q\dfrac{1}{k_i+1} < \dfrac{q-3n+2}{2q-5n+10}\biggl(\dfrac{2(2q+n-3)}{3n(n+1)}-3\biggl)$,
then there are two maps $f^s,f^t\ (1\le s<t\le 3)$ and $n+1$ indices $1\le i_1<i_2<\cdots <i_{n+1}\le q$ such that
$$\dfrac{(f^s,\widetilde{a}_{i_1}^1)}{(f^t,\widetilde{a}_{i_1}^2)} =\cdots = \dfrac{(f^s,\widetilde{a}_{i_{n+1}}^1)}{(f^t,\widetilde{a}_{i_{n+1}}^2)}.$$
\end{theorem}

\section{Basic notions and auxiliary results from Nevanlinna theory}

\noindent
{\bf (a)} Counting function of divisor.

For $z = (z_1,\dots,z_m) \in \C^m$, we set
$\Vert z \Vert = \Big(\sum\limits_{j=1}^m |z_j|^2\Big)^{1/2}$
and define
\begin{align*}
B(r) &= \{ z \in \C^m ; \Vert z \Vert < r\},\quad
S(r) = \{ z \in \C^m ; \Vert z \Vert = r\},\\
d^c & = \dfrac{\sqrt{-1}}{4\pi}(\overline \partial - \partial),\quad
\sigma = \big(dd^c \Vert z\Vert^2\big)^{m-1},\\
\eta &= d^c \text{log}\Vert z\Vert^2 \land
\big(dd^c\text{log}\Vert z \Vert\big)^{m-1}.
\end{align*}

Thoughout this paper, we denote by $\mathcal M$ the set of all meromorphic functions on $\C^m$. A divisor $E$ on $\C^m$ is given by a formal sum $E=\sum\mu_{\nu}X_{\nu}$,
where $\{X_\nu\}$ is a locally family of distinct irreducible analytic hypersurfaces in $\C^m$ and $\mu_{\nu}\in \Z$. We define the support of the divisor
$E$ by setting $\supp (E)=\cup_{\nu\ne 0} X_\nu$.
Sometimes, we identify the divisor $E$ with a function $E(z)$ from $\C^m$
into $\Z$ defined by $E(z):=\sum_{X_{\nu}\ni z}\mu_\nu$.

Let $M,k$ be a positive integer or $+\infty$. We define the truncated divisor $E^{[M]}$ and $E_{\le k}^{[M]}$ by
\begin{align*}
E^{[M]}:= \sum_{\nu}\min\{\mu_\nu, M \}X_\nu ,\\
 E_{\leq k}^M:=
\begin{cases}
0,  & \text  { if } E(z) > k,\\
E^M,  & \text  { if } E(z) \leq k.
\end{cases}
\end{align*}
and the {\it truncated counting function to level $M$} of $E$ by
\begin{align*}
N^{[M]}(r,E) := \int\limits_1^r \frac{n^{[M]}(t,E)}{t^{2m-1}}dt\quad
(1 < r < +\infty),
\end{align*}
Similarly, we define $N(r,E_{>k}^{[M]})$ and $N(r,E_{\leq k}^{[M]})$ and denote them by $N_{>k}^{[M]}(r,E)$ and $N_{\leq k}^{[M]}(r,E)$ respectively.
where
\begin{align*}
n^{[M]}(t,E): =
\begin{cases}
\int\limits_{\supp (E) \cap B(t)} E^{[M]}\sigma &\text{ if } m \geq 2,\\
\sum_{|z| \le t} E^{[M]}(z)&\text{ if } m = 1.
\end{cases}
\end{align*}
Similarly, we define $n_{>k}^{[M]}(t,E)$ and $n_{\leq k}^{[M]}(t,E)$.

We omit the character $^{[M]}$ if $M=+\infty$.

 For an analytic hypersurface $E$ of $\C^m$, we may consider it as a reduced divisor and denote by $N(r,E)$ its counting function.

Let $\varphi$ be a nonzero meromorphic function on $\C^m$.
We denote by $\nu^0_{\varphi}$ (resp. $\nu^{\infty}_{\varphi}$) the divisor of zeros (resp. divisor of poles) of $\varphi$. The divisor of $\varphi$ is defined by
$$\nu_{\varphi}=\nu^0_{\varphi}-\nu^{\infty}_{\varphi}.$$

We have the following Jensen's formula:
\begin{align*}
N(r,\nu^0_{\varphi}) - N(r,\nu^{\infty}_{\varphi}) =
\int\limits_{S(r)} \text{log}|\varphi| \eta
- \int\limits_{S(1)} \text{log}|\varphi| \eta .
\end{align*}

For convenience, we will write $N_{\varphi}(r)$ and $N^{[M]}_{\varphi}(r)$ for $N(r,\nu^0_{\varphi})$ and $N^{[M]}(r,\nu^0_{\varphi})$, respectively.

\noindent
{\bf (b)} The first main theorem.

Let $f$ be a meromorphic mapping of $\C^m$ into $\P^n(\C)$. For arbitrary fixed homogeneous coordinates $(w_0: \cdots : w_n)$ of $\P^n(\C )$, we take a reduced representation $f = (f_0 : \cdots : f_n)$, which means that each $f_i$ is holomorphic function on $\C^m$ and $f(z) = (f_0(z) : \cdots : f_n(z))$ outside the analytic set $I(f):=\{ z ; f_0(z) = \cdots = f_n(z) = 0\}$ of codimension at least $2$.

Denote by $\Omega$ the Fubini Study form of $\P^n(\C )$. The characteristic function of $f$ (with respect to $\Omega$) is defined by
\begin{align*}
T_f(r) := \int_1^r\dfrac{dt}{t^{2m-1}}\int_{B(t)}f^*\Omega\wedge\sigma ,\quad\quad 1 < r < +\infty.
\end{align*}
By Jensen's formula we have
\begin{align*}
T_f(r)=\int_{S(r)}\log ||f||\eta +O(1),
\end{align*}
where $\Vert f  \Vert = \max \{ |f_0|,\dots,|f_n|\}$.

Let $a$ be a meromorphic mapping of $\C^m$ into $\P^n(\C)^*$ with reduced representation
$a = (a_0 : \dots : a_n)$. We define
$$m_{f,a}(r)=\int\limits_{S(r)} \text{log}\dfrac {||f||\cdot ||a||}{|(f,a)|}\eta -
\int\limits_{S(1)}\text{log}\dfrac {||f||\cdot ||a||}{|(f,a)|}\eta,$$
where $\Vert a \Vert = \big(|a_0|^2 + \dots + |a_n|^2\big)^{1/2}$ and $(f,a)=\sum_{i=0}^nf_i\cdot a_i.$

Let $f$ and $a$ be as above. If $(f,a)\not \equiv 0$, then the first main theorem for moving hyperplaness in value distribution theory states
$$T_f(r)+T_a(r)=m_{f,a}(r)+N_{(f,a)}(r)+O(1)\ (r>1).$$

For a meromorphic function $\varphi$ on $\C^m$, the proximity function $m(r,\varphi)$ is defined by
$$ m(r,\varphi) = \int\limits_{S(r)} \log^+ |\varphi| \eta , $$
where $\log^+ x = \max \big\{ \log x, 0\big\}$ for $x \geqslant 0$.
The Nevanlinna's characteristic function is defined by
$$T(r, \varphi ) = N(r, \nu^{\infty}_\varphi) + m(r,\varphi ).$$
We regard $\varphi$ as a meromorphic mapping of $\C^m$ into $\P^1(\C )^*$, there is a fact that
$$ T_\varphi (r)=T(r,\varphi )+O (1). $$
{\bf (c)} Lemma on logarithmic derivative.

 As usual, by the notation $``|| \ P"$  we mean the assertion $P$ holds for all $r \in [0,\infty)$ excluding a Borel subset $E$ of the interval $[0,\infty)$ with $\int_E dr<\infty$. Denote by $\Z_+$ the set of all nonnegative integers. The lemma on logarithmic derivative in Nevanlinna theory is stated as follows.
\begin{lem}[{see \cite[Lemma 3.11]{Shi}}]\label{2.1}
Let $f$ be a nonzero meromorphic function on $\C^m.$ Then
$$\biggl|\biggl|\quad m\biggl(r,\dfrac{\mathcal{D}^\alpha (f)}{f}\biggl)=O(\log^+T_f(r))\ (\alpha\in \Z^m_+).$$
\end{lem}

\noindent
{\bf (d)} Family of moving hyperplanes.

We assume that thoughout this paper, the homogeneous coordinates of $\P^n(\C)$ is chosen so that for each given meromorphic mapping $a=(a_0:\cdots :a_n)$ of $\C^m$ into $\P^n(\C)^*$ then $a_{0}\not\equiv 0$. We set
$$ \tilde a_i=\dfrac{a_i}{a_0}\text{ and }\tilde a=(\tilde a_0:\tilde a_1:\cdots:\tilde a_n).$$
Let $f:\C^m\rightarrow\P^n(\C)$ be a meromorphic mapping with the reduced representation $f=(f_0:\cdots :f_n).$  We put $(f,a):=\sum_{i=0}^{n}f_ia_{i}$ and $(f,\tilde a):=\sum_{i=0}^{n}f_i\tilde a_{i}.$

Let $\{a_i\}_{i=1}^q$ be $q$ meromorphic mappings of $\C^m$ into $\P^n(\C)^*$  with reduced representations $a_i=(a_{i0}:\cdots :a_{in})\ (1\le i\le q).$ We denote by  $\mathcal R(\{a_i\})$ (for brevity we will write $\mathcal R$ if there is no confusion) the smallest subfield of $\mathcal M$ which contains $\C$ and all ${a_{i_j}}/{a_{i_k}}$ with $a_{i_k}\not\equiv 0.$

\begin{definition}\label{2.2}
The family $\{a_i\}_{i=1}^q$ is said to be in general position  if $\dim (\{a_{i_0},\ldots ,a_{i_n}\})_{\mathcal {M}}=n+1$ for any $1\le i_0\le\cdots\le i_n\le q$, where $(\{a_{i_0},\ldots ,a_{i_n}\})_{\mathcal {M}}$ is the linear span of $\{a_{i_0},\ldots ,a_{i_N}\}$ over the field $\mathcal{ M}.$
\end{definition}

\begin{theorem}[{The second main theorem \cite[Corollary 1.2]{Q}}]\label{2.3}
Let $f :\C^m \to \P^n(\C)$ be a meromorphic mapping. Let $\{a_i\}_{i=1}^q \ (q\ge 2n+1)$ be  meromorphic mappings of $\C^m$ into $\P^n(\C)^*$ in general position such that $(f,a_i)\not\equiv 0\ (1\le i\le q).$

$\mathrm{(a)}$ If $q\ge 3n+3$ then
$$|| \dfrac{2q}{3(n+1)}T_f(r) \le \sum_{i=1}^q N_{(f,a_i)}^{[n]}(r) + o(T_f(r)) + O(\max_{1\le i \le q}T_{a_i}(r)).$$

$\mathrm{(b)}$ If $q< 3n+3$ then
$$|| \dfrac{q-n+1}{n+2}T_f(r) \le \sum_{i=1}^q N_{(f,a_i)}^{[n]}(r) + o(T_f(r)) + O(\max_{1\le i \le q}T_{a_i}(r)).$$
\end{theorem}

\section{Proof of Theorem \ref{1.1}}
Assume that $\sum_{v=1}^q\dfrac{1}{k_v+1}< \biggl(\dfrac{2q}{3n(n + 1)}-\dfrac{2q}{q+2n-2}\biggl).$ Suppose that  the conclussion \ref{new2} does not hold.  By changing indices if necessary, we may assume that
$$\underbrace{\dfrac{(f^1,\widetilde{a}_1^1)}{(f^2,\widetilde{a}_1^2)}\equiv \dfrac{(f^1,\widetilde{a}_2^1)}{(f^2,\widetilde{a}_2^2)}\equiv \cdots\equiv \dfrac{(f^1,\widetilde{a}_{v_1}^1)}
{(f^2,\widetilde{a}_{v_1}^2)}}_{\text { group } 1}\not\equiv
\underbrace{\dfrac{(f^1,\widetilde{a}_{v_1+1}^1)}{(f^2,\widetilde{a}_{v_1+1}^2)}\equiv \cdots\equiv \dfrac{(f^1,\widetilde{a}_{v_2}^1)}{(f^2,\widetilde{a}_{v_2}^2)}}_{\text { group } 2}$$
$$\not\equiv \underbrace{\dfrac{(f^1,\widetilde{a}_{v_2+1}^1)}{(f^2,\widetilde{a}_{v_2+1}^2)}\equiv \cdots\equiv \dfrac{(f^1,\widetilde{a}_{v_3}^1)}{(f^2,\widetilde{a}_{v_3}^2)}}_{\text { group } 3}\not\equiv \cdots\not\equiv \underbrace{\dfrac{(f^1,\widetilde{a}_{v_{s-1}+1}^1)}{(f^2,\widetilde{a}_{v_{s-1}+1}^2)}\equiv\cdots \equiv
\dfrac{(f^1,\widetilde{a}_{v_s}^1)}{(f^2,\widetilde{a}_{v_s}^2)}}_{\text { group } s},$$
where $v_s=q.$

For each $1\le i \le q,$ we set
\begin{equation*}
\sigma (i)=
\begin{cases}
i+n& \text{ if $i+n\leq q$},\\
i+n-q&\text{ if  $i+n> q$}
\end{cases}
\end{equation*}
and
$$P_i=(f^1,\widetilde{a}_i^1)(f^2,\widetilde{a}_{\sigma (i)}^2)-(f^2,\widetilde{a}_i^2)(f^1,\widetilde{a}_{\sigma (i)}^1).$$
By supposition, the number of elements of each group is at most $n$. Hence $\dfrac{(f^1,\widetilde{a}_i^1)}{(f^2,\widetilde{a}_i^2)}$ and $\dfrac{(f^1,\widetilde{a}_{\sigma (i)}^1)}{(f^2,\widetilde{a}_{\sigma (i)}^2)}$ belong to distinct groups. This means that $P_i\not\equiv 0\ (1\le i\le q)$.

Fix an index $i$ with $1\le i \le q.$ It is easy to see that
\begin{align*}
\nu_{P_i}(z) \ge \min\{\nu_{(f^1,\widetilde{a}_i^1)},\nu_{(f^2,\widetilde{a}_i^2)}\}+\min\{\nu_{(f^1,\widetilde{a}_{\sigma (i)}^1)},\nu_{(f^2,\widetilde{a}_{\sigma (i)}^2)}\}
+\sum_{\underset{v\ne i,\sigma (i)}{v=1}}^{q}\nu_{(f^1,\widetilde{a}_v^1)}^{[1]}(z)
\end{align*}
outside a finite union of analytic sets of dimension $\le m-2.$ Since $\min\{a,b\}+n\ge\min\{a,n\}+\min\{b,n\}$ for all positive integers $a$ and $b$,
the above inequality implies that
\begin{align*}
N_{P_i}(r)\geq \sum_{v=i,\sigma (i)}\left (
N^{[n]}_{(f^1,\widetilde{a}_v^1),\le k_v}(r)+N^{[n]}_{(f^2,\widetilde{a}_v^2),\le k_v}(r)-nN^{[1]}_{(f^1,\widetilde{a}_v^1),\le k_v}(r) \right )
+\sum_{\underset{v\ne i,\sigma (i)}{v=1}}^{q}N^{[1]}_{(f^1,\widetilde{a}_v^1),\le k_v}(r).
\end{align*}
On the other hand, by the Jensen formula, we have
\begin{align*}
N_{P_i}(r)=&\int_{S(r)}\log |P_i|\eta + O(1)\\
\le &\int_{S(r)}\log (|(f^1,\widetilde{a}_i^1)|^2+|(f^1,\widetilde{a}_{\sigma (i)}^1|^2)^{\frac{1}{2}}\eta
+ \int_{S(r)}\log (|(f^2,\widetilde{a}_i^2)|^2+|(f^2,\widetilde{a}_{\sigma (i)}^2|^2)^{\frac{1}{2}}\eta +O(1)\\
\le &T_{f^1}(r)+T_{f^2}(r) +o(T_{f^1}(r)+T_{f^2}(r)).
\end{align*}
This implies that
\begin{align*}
T_{f^1}(r)+T_{f^2}(r)\ge &\sum_{v=i,\sigma (i)}\left (
N^{[n]}_{(f^1,\widetilde{a}_v^1),\le k_v}(r)+N^{[n]}_{(f^2,\widetilde{a}_v^2),\le k_v}(r)-nN^{[1]}_{(f^1,\widetilde{a}_v^1),\le k_v}(r) \right )\\
+ &\sum_{\underset{v\ne i,\sigma (i)}{v=1}}^{q}N^{[1]}_{(f^1,\widetilde{a}_v^1),\le k_v}(r)
+o(T_{f^1}(r)+T_{f^2}(r)).
\end{align*}
Summing-up both sides of the above inequality over $i=1,\ldots ,q$, we have
\begin{align*}
q(T_{f^1}(r)+T_{f^2}(r))\ge &2\sum_{v=1}^q\left (
N^{[n]}_{(f^1,\widetilde{a}_v^1),\le k_v}(r)+N^{[n]}_{(f^2,\widetilde{a}_v^2),\le k_v}(r) \right )\\
&+ (q-2n-2)\sum_{v=1}^{q}N^{[1]}_{(f^1,\widetilde{a}_v^1),\le k_v}(r)+o(T_{f^1}(r)+T_{f^2}(r))\\
\ge & (2+\frac{q-2n-2}{2n})\sum_{v=1}^q\left (
N^{[n]}_{(f^1,\widetilde{a}_v^1),\le k_v}(r)+N^{[n]}_{(f^2,\widetilde{a}_v^2),\le k_v}(r) \right )\\
&+o(T_{f^1}(r)+T_{f^2}(r))\\
\end{align*}
We get
\begin{align*}
\dfrac{2qn}{q+2n-2}(T_{f^1}(r)+T_{f^2}(r))&\ge \sum_{v=1}^q\left (
N^{[n]}_{(f^1,\widetilde{a}_v^1),\le k_v}(r)+N^{[n]}_{(f^2,\widetilde{a}_v^2),\le k_v}(r) \right )+o(T_{f^1}(r)+T_{f^2}(r))\\
&= \sum_{v=1}^q(N^{[n]}_{(f^1,\widetilde{a}_v^1)}(r)+N_{(f^2,\widetilde{a}_v^2)}^{[n]}(r)
- N_{(f^1,\widetilde{a}_v^1),> k_v}^{[n]}(r)-N_{(f^2,\widetilde{a}_v^2),> k_v}^{[n]}(r))\\
&+o(T_{f^1}(r)+T_{f^2}(r))\\
&\ge \sum_{v=1}^q(N^{[n]}_{(f^1,a_v^1)}(r)+N_{(f^2,a_v^2)}^{[n]}(r)
- N_{(f^1,\widetilde{a}_v^1),> k_v}^{[n]}(r)-N_{(f^2,\widetilde{a}_v^2),> k_v}^{[n]}(r))\\
&+o(T_{f^1}(r)+T_{f^2}(r))\\
\end{align*}
By theorem \ref{2.3}, we have
\begin{align*}
\biggl|\biggl|  \dfrac{2q}{3(n+1)}(T_{f^1}(r)+T_{f^2}(r))&\le \sum_{v=1}^q(N^{[n]}_{(f^1,a_v^1)}(r)+N_{(f^2,a_v^2)}^{[n]}(r))+o(T_{f^1}(r)+T_{f^2}(r))\\
 \end{align*}
From the above inequalities, we have
\begin{align*}
\biggl(\frac{2q}{3(n + 1)}&-\frac{2qn}{q+2n-2}\biggl)(T_{f^1}(r)+T_{f^2}(r)) \\
&\leq \sum_{v=1}^q \biggl (N_{(f^1,\widetilde{a}_v^1),> k_v}^{[n]}(r)+N_{(f^2,\widetilde{a}_v^2),> k_v}^{[n]}(r)\biggl)\\
&\ + o(T_{f^1}(r)+T_{f^2}(r))\\
&\leq  \sum_{v=1}^q\dfrac{n}{k_v+1}(N_{(f^1,\widetilde{a}_v^1)}(r)+N_{(f^2,\widetilde{a}_v^2)}(r)) \\
&\ + o(T_{f^1}(r)+T_{f^2}(r))\\
&\leq n\sum_{v=1}^q\dfrac{1}{k_v+1}(T_{f^1}(r)+T_{f^2}(r))+ o(T_{f^1}(r)+T_{f^2}(r))
\end{align*}
Letting $r \to \infty$, we get
$$\biggl(\dfrac{2q}{3n(n + 1)}-\dfrac{2q}{q+2n-2}\biggl) \leq \sum_{v=1}^q\dfrac{1}{k_v+1}.$$
This is a contradiction.

Then the supposition is impossible. Hence the theorem is proved. \hfill $\square$
\section{Proof of Theorem \ref{1.2}}

In order to prove Theorem \ref{1.2}, we need the following.

\noindent
{\bf 3.1.} Let $f^1,f^2,f^3:\C^m \to \P^n(\C)$ be three meromorphic mappings.  Let $k_i$ $(1\le i\le q)$ be positive integers or $\infty$. Let $\{a_i^t\}_{i=1}^{q}$ $(t=1,2,3)$ be 3 families of  moving hyperplanes in $\P^n(\C)$ in general position such that $a_i^t$ be "slowly'' with respect to $f^t$ and $\dim\ \{z\in\C^m: \nu_{(f^t,a^t_i),\le k_i}.\nu_{(f^t,a^t_j),\le k_j}>0\} \leq m-2\quad (1\le i<j\le q, 1\le t \le 3)$, we put
$$T(r)=\sum_{t=1}^{3} T_{f^t}(r).$$
Assume that $a_i^t$ has a reduced representation $a_i^t=(a_{i0}^t:\cdots :a_{in}^t).$ By changing the homogeneous coordinate system of $\P^n(\C),$
we may assume that $a_{i0}^t\not \equiv 0\ (1\le i \le q, 1\le t \le 3).$ For each $c=(c_1,...,c_q) \in \C^q \setminus \{0\}$, we set
\begin{align*}
a^t_c:=(\sum_{i=1}^qc_i\widetilde{a}^t_{i0},...,\sum_{i=1}^qc_i\widetilde{a}^t_{in}), \text{ } ||a^t_c||:=(\sum_{j=0}^n|\sum_{i=1}^qc_i\widetilde{a}^t_{ij}|^2)^{\frac{1}{2}}\\
(f^t,a^t_c):=\sum_{j=0}^n\sum_{i=1}^qc_i\widetilde{a}_{ij}f^t_j=\sum_{i=1}^qc_i(f^t,\widetilde{a}_i)\quad (1\le t\le 3)
\end{align*}
We denote by $\beta$ the union of all irreducible components with dimension $m-1$ of the analytic set $\bigcap_{i=1}^q Zero(f^t,a^t_i)$ $(1\le t\le 3)$. Then $\beta$ is either an analytic set of pure dimension $m-1$ or empty set. With $c \in \C^q,$ we denote by $S_c^{jt}$ the closure of set $(Zero(f^t,a^t_j)\cap Zero(f^t,a^t_c))\setminus \beta$. Then $S^{jt}_c$ is an analytic set. We also denote by $\mathcal{C}$ the set of all $c \in \C^q\setminus\{0\}$ such that $\dim \mathcal{S}_c^{jk}\le m-2$
\begin{lemma}
$\mathcal {C}$ is dense in $\C^{q}.$
\end{lemma}
\textbf{\textit{Proof.}} \ For $1\le i\le q, 1\le t\le 3$ and for each irreducible component $\nu$ of the analytic set $\zero (f^{t},a^t_i)$ with $\nu\not\subset\beta$, we set
$$
V^{it}_{\nu}=\{c=(c_1,\ldots,c_q)\in\C^{q}\ :\ (f^t,a^t_c)(z)=0,\ \forall z\in\nu\}.
$$
Then, $V^{it}_{\nu}$ is an complex vector subspace of $\C^q$. Since $\nu\not\subset\beta$, there exists an index $j$ such that $\nu\not\subset \zero (f^{t},a^t_j)$. Therefore the element $c=(0,\ldots,0,\underset{j-th}{1},0,\ldots,0)$ does not belong to $V^{it}_{\nu}$. Hence $\dim V^{it}_{\nu}\le q-1$. Let $K=\bigcup_{i=1}^{q}\bigcup_{t=1}^{3}\bigcup_{\nu}V^{it}_{\nu}$. Then $K$ is a union of at most a countable number of $(q - 1)$-dimensional complex vector subspaces in $\C^q$. It is easy to see that $\mathcal C\supset \C^q\setminus K$. Therefore $\mathcal C$ is dense in $\C^q$. The lemma is proved.\hfill$\square$
\begin{lemma}\label{4.0}
For every $c \in \mathcal{C}$, we put $F^{jt}_{c}:=\dfrac{(f^t,\widetilde{a}_j^t)}{(f^t,a_c^t)}\quad (1 \le j \le q, \ 1\le t\le 3).$ Then
$$||T(r,F_c^{jt})\le T_{f^t}(r)+o(T(r))$$
\end{lemma}\label{4.1}
\textbf{\textit{Proof.}}\ Let $h$ be a meromorphic function on $\C^m$ such that $\big (h(f^t, \tilde a^t_j):h(f^t, a^t_c)\big )$ is a reduced representation of a meromorphic mapping into $\P^{1}(\C )$. It is easy to see that
$$
\nu^0_h\le \sum_{i=1}^q\nu_{a_{j0}}.
$$
This implies that
$$
||\ N_h(r)\le \sum_{j=1}^qN_{a^t_{j0}}(r)\le\sum_{j=1}^qT_{a^t_j}(r)=o(T(r)).
$$
By the definition of the characteristic function and by Jensen formula, we have
\begin{align*}
||\ T(r, F^{jt}_c)&=\int\limits_{S(r)}\log\left (|h(f^t, \tilde a^t_j)|^2+|h(f^t, a^t_c)|^2\right )^{\frac{1}{2}}\eta \\
&\le \int\limits_{S(r)}\log ||f^t||\eta +\int\limits_{S(r)}\log |h|\eta +\int\limits_{S(r)}\log (||\tilde a^t_j||^2+||a^t_c||^2)^{\frac{1}{2}}\eta +O(1)\\
&\le T_{f^t}(r)+N_h(r)+\int\limits_{S(r)}\log^+ ||\tilde a^t_j||\eta +\int\limits_{S(r)}\log^+ ||\tilde a^t_c||\eta +O(1)\\
&\le T_{f^t}(r)+\sum_{i=1}^{n}\int\limits_{S(r)}\log^+ |\dfrac{a^t_{ji}}{a_{j0}}|\eta +\sum_{v=1}^{q}\sum_{i=1}^{n}\int\limits_{S(r)}\log^+ |\dfrac{a^t_{vi}}{a_{v0}}|\eta+o(T(r))\\
&=  T_{f^t}(r)+\sum_{i=1}^{n}m(r,\dfrac{a^t_{ji}}{a^t_{j0}})+\sum_{v=1}^{q}\sum_{i=1}^{n}m(r,\dfrac{a^t_{vi}}{a_{v0}})+o(T(r))\\
&\le T_{f^t}(r)+nT_{a^t_j}(r)+n\sum_{v=1}^{q}T_{a^t_v}(r)+o(T(r))= T_{f^t}(r)+o(T(r)).
\end{align*}
The lemma is proved.\hfill$\square$

\begin{definition}[{see \cite[p. 138]{Fu3}}]
Let $F_1,F_2 ,F_3$ be nonzero meromorphic functions on $\C^m$. Take a set $\alpha=(\alpha_1,...,\alpha_m) \in (\Z^+)^m$ with $|\alpha|=\sum_{i=1}^m\alpha_i =1$.  We define Cartan's auxiliary function by
$$\Phi^\alpha \equiv \Phi^\alpha(F_1,F_2,F_3):=F_1F_2F_3\left | \begin {array}{cccc}
1&1&1\\
\frac {1}{F_1}&\frac {1}{F_2} &\frac {1}{F_3}\\
\mathcal {D}^{\alpha}(\frac {1}{F_1}) &\mathcal {D}^{\alpha}(\frac {1}{F_2}) &\mathcal {D}^{\alpha}(\frac {1}{F_3})\\
\end {array}
\right|
$$
\end{definition}
By simple computation, we have
\begin{align}\label{new1}
\Phi^\alpha(F_1,F_2,F_3)=F_1\biggl (\dfrac{\mathcal D^\alpha F_2}{F_2}-\dfrac{\mathcal D^\alpha F_3}{F_3}\biggl)+F_2\biggl (\dfrac{\mathcal D^\alpha F_3}{F_3}-\dfrac{\mathcal D^\alpha F_1}{F_1}\biggl)+F_3\biggl (\dfrac{\mathcal D^\alpha F_1}{F_1}-\dfrac{\mathcal D^\alpha F_2}{F_2}\biggl).
\end{align}
\begin{lemma}[{see \cite[Proposition 3.4]{Fu3}}]\label{4.4}
 If $\Phi^\alpha(F,G,H)=0$ and $\Phi^\alpha(\frac {1}{F},\frac {1}{G},\frac {1}{H})=0$ for all $\alpha$ with $|\alpha|=1$, then one of the following assertions holds :

(i) \ $F=G, G=H$ or $H=F$

(ii) \ $\frac {F}{G},\frac {G}{H}$ and $\frac {H}{F}$ are all constant.
\end{lemma}
%\begin{lemma}[{see \cite[Proposition 3.3]{Fu3}}]\label{4.5}
%Suppose that $\Phi^\alpha(F_1,...,F_{M+1}) \not \equiv 0$ with $|\alpha|\le \frac{M(M+1)}{2}$. If
%$$\nu^{[d]}:=\min(\nu_{F_1},d)=\min(\nu_{F_2},d)=\min(\nu_{F_{M+1}},d)$$
%for some $d \ge |\alpha|$, then $\nu_{\Phi^\alpha}(z_0) \ge \min(\nu^{[d]}(z_0),d-|\alpha|)$ for every $z_0 \in F_1^{-1}(0)\backslash A$, where A is an analytic subset of $codim A \ge 2$.
%\end{lemma}\begin{lemma}[{see \cite[Proposition 3.7]{Fu3}}]\label{4.6}
%Let $F_1,...,F_{M+1}$ and $\alpha$ satisfy the same assumptions as in Lemma $\ref{4.5}$, if $F_1=...=F_{M+1} \not \equiv 0, \infty$ on an analytic subset A of pure dimension $n-1$, then $\nu_{\Phi^\alpha}(z_0) \ge M$ for every $z_0 \in A$.
%\end{lemma}
\begin{lemma}[{see \cite[Lemma 4.7]{TQ05}}]\label{4.7}
Suppose that there exists $\Phi^\alpha=\Phi^\alpha(F_{c}^{j_01},F_{c}^{j_02},F_{c}^{j_3})\not\equiv 0$ for some $c \in \mathcal{C},$ $|\alpha|=1$. Then, for each $1 \le t \le 3$, the following holds:
\begin{align*}
2\sum_{j=1}^3N^{[1]}_{(f^1,\widetilde{a}^1_{j}),\le k_{j}}(r)&+\sum_{t=1}^3N^{[n]}_{(f^t,\widetilde{a}^t_{j_0}),\le k_{j_0}}(r)
-(2n+3)N^{[1]}_{(f^1,\widetilde{a}^1_{j_0}),\le k_{j_0}}(r)-2\sum_{t=1}^3N^{[1]}_{(f^t,\widetilde{a}^t_{j_0}),>k_{j_0}}(r)\\
&\le N_{\Phi^\alpha}(r)+ o(T(r))\le T(r)+\sum_{t=1}^3N^{[1]}_{(f^t,\widetilde{a}^t_{j_0}),>k_{j_0}}(r)+o(T(r)).
\end{align*}
\end{lemma}
{\bf Proof.}
(a)Firstly, we will prove the first inequality. We set
\begin{align*}
&\mathcal{A}=\{z \in \C^m: \nu^0_{(f^t,\widetilde{a}^t_{j_0})}>0\}\\
&\mathcal{V}=\{z \in \C^m: \nu^0_{(f^t,\widetilde{a}^t_i),\le k_i}.\nu^0_{(f^t,\widetilde{a}^t_j),\le k_j}>0\}\quad (1\le i<j\le q)
\end{align*}
Then $\mathcal{V}$ is an analytic set of codimension at least 2. We also set
$$
\mathcal{D}=\bigcup_{i=1}^q\{z \in \C^m: \nu^0_{(f^t,\widetilde{a}^t_i),\le k_i}>0\} \text{ and } \mathcal{S}=(\bigcup_{i=1}^q\bigcup_{t=1}^3S_c^{it})\cup \mathcal{A}\cup \beta
$$
Let $z_0$ be a regular point of the analytic set $\mathcal{D}$ such that $z_0 \not \in \mathcal{V}\cup \mathcal{S}$. There are three cases:

{\bf Case 1.} $z_0 \not \in \mathcal{A}$. Let $\nu$ be the irreducible component of $\mathcal D$ which contains $z_0$. Then, there exist a neighborhood $U$ of $z_0$ and a holomorphic function $h$ on $U$ such that $dh$ has nonzero point and $U\cap\zero h=\nu$. Moreover, we may assume that $U\cap (\mathcal V\cup\mathcal S\cup\mathcal{A})=\emptyset$. Since $\dfrac{(f^t,\widetilde{a}^t_i)}{(f^s,\widetilde{a}^s_i)}=\dfrac{(f^t,\widetilde{a}^t_j)}{(f^s,\widetilde{a}^s_j)}$ for all $z\in\nu$, $1\le i \ne j \le q$, $1\le t \ne s \le 3$, there exist holomorphic functions $\varphi_v$ defined on $U$ such that $F^{cv}_{j_0}=h\varphi_v$ on U $(1\le v \le 3)$

Then, we rewrite the function $\Phi^\alpha$ on $U$ as follows
\begin{align*}
\Phi^\alpha(F_{c}^{j_01},F_{c}^{j_02},F_{c}^{j_03}):&=F_{c}^{j_01}F_{c}^{j_02}F_{c}^{j_03}\left | \begin {array}{cccc}
1&1&1\\
F_{j_0}^{c1}&F_{j_0}^{c2} &F_{j_0}^3\\
\mathcal {D}^{\alpha}(F_{j_0}^{c1}) &\mathcal {D}^{\alpha}(F_{j_0}^{c2})  &\mathcal {D}^{\alpha}(F_{j_0}^{c3})\\
\end {array}
\right|\\
&=F_{c}^{j_01}F_{c}^{j_02}F_{c}^{j_03}\left | \begin {array}{cccc}
F^{c2}_{j_0}-F^{c1}_{j_0} &F^{c3}_{j_0}-F^{c1}_{j_0}\\
\mathcal {D}^{\alpha}(F^{c2}_{j_0}-F^{c1}_{j_0}) &\mathcal {D}^{\alpha}(F^{c3}_{j_0}-F^{c1}_{j_0})\\
\end {array}
\right|\\
&=F_{c}^{j_01}F_{c}^{j_02}F_{c}^{j_03}h^2\left | \begin {array}{cccc}
\varphi_2-\varphi_1 &\varphi_{3}-\varphi_1\\
\mathcal {D}^{\alpha}(\varphi_2-\varphi_1) &\mathcal {D}^{\alpha}(\varphi_3-\varphi_1)\\
\end {array}
\right|\\
\end{align*}
\begin{align*}
 \nu^0_{\Phi^\alpha}(z_0)\ge 2&\ge 2\sum_{j=1}^q\min\{1,\nu_{(f^1,\widetilde{a}^1_j),\le k_j}(z_0)\}+\sum_{t=1}^3\min\{n,\nu_{(f^t,\widetilde{a}^t_{j_0}),\le k_{j_0}}(z_0)\}\\
&\ -(2n+3)\min\{1,\nu_{(f^1,\widetilde{a}^1_{j_0}),\le k_{j_0}}(z_0)\}-2\sum_{t=1}^3\min\{1,\nu_{(f^t,\widetilde{a}^t_{j_0}),>k_{j_0}}(z_0)\}
\end{align*}
%{\bf Case 2.} $z_0 \in Zero(f^t,a^t_{j_0})$.
%So that $|| N_{\Phi^\alpha}(r)\ge N^{[d-|\alpha|]}_{(f^t,a^t_{j_0}),\le k_{j_0}}(r)+M\sum_{j\ne j_0}N^{[1]}_{(f^t,a^t_{j_0}),\le k_{j_0}}(r)$ $(1\le t \le M+1)$
  {\bf Case 2.} $z_0 \in \{z \in \C^m: \nu^0_{(f^t,\widetilde{a}^t_{j_0}),\le k_{j_0}}>0\}$. Without loss of generality, we may assume that $0<\nu_{(f^1,\widetilde{a}^1_{j_0}),\le k_{j_0}}(z_0)\le \nu_{(f^2,\widetilde{a}^2_{j_0}),\le k_{j_0}}(z_0)\le \nu_{(f^3,\widetilde{a}^3_{j_0}),\le k_{j_0}}(z_0)$.
  \begin{align*}
\Phi^\alpha&=F_{c}^{j_01}\left | \begin {array}{cccc}
F^{c2}_{j_0}(F^{c2}_{j_0}-F^{c1}_{j_0}) &F^{c3}_{j_0}(F^{c3}_{j_0}-F^{c1}_{j_0})\\
F^{c2}_{j_0}\mathcal {D}^{\alpha}(F^{c2}_{j_0}-F^{c1}_{j_0}) &F^{c3}_{j_0}\mathcal {D}^{\alpha}(F^{c3}_{j_0}-F^{c1}_{j_0})\\
\end {array}
\right|\\
&=F_{c}^{j_01}\biggl(F^{c2}_{j_0}(F^{c2}_{j_0}-F^{c1}_{j_0}).F^{c3}_{j_0}\mathcal {D}^{\alpha}(F^{c3}_{j_0}-F^{c1}_{j_0})-F^{c3}_{j_0}(F^{c3}_{j_0}-F^{c1}_{j_0}).F^{c2}_{j_0}\mathcal {D}^{\alpha}(F^{c2}_{j_0}-F^{c1}_{j_0})\biggl) \\
\end{align*}
Because of the assumption, we see that $F^{c2}_{j_0}(F^{c2}_{j_0}-F^{c1}_{j_0})$ and $F^{c3}_{j_0}(F^{c3}_{j_0}-F^{c1}_{j_0})$ are holomorphic on a neighborhood of $z_0$. Moreover, we have
\begin{align*}
&\nu^\infty_{F^{c2}_{j_0}\mathcal{D}^{\alpha}(F^{c2}_{j_0}-F^{c1}_{j_0})}(z_0) \le |\alpha|=1\\
&\nu^\infty_{F^{c3}_{j_0}\mathcal{D}^{\alpha}(F^{c3}_{j_0}-F^{c1}_{j_0})}(z_0) \le |\alpha|=1
\end{align*}
 Therefore
\begin{align*}
\nu^0_{\Phi^\alpha}(z_0)&\ge \nu_{F_c^{j_01}}(z_0)-1=\min_{1\le t\le 3} \nu_{(f^t,\widetilde{a}^t_{j_0}),\le k_{j_0}}(z_0)-\min\{1,\nu_{(f^1,\widetilde{a}^1_{j_0}),\le k_{j_0}}(z_0)\}\\
&\ge\sum_{t=1}^3\min\{n,\nu_{(f^t,\widetilde{a}^t_{j_0}),\le k_{j_0}}(z_0)\}-2n\min\{1,\nu_{(f^1,\widetilde{a}^1_{j_0}),\le k_{j_0}}(z_0)\}-\min\{1,\nu_{(f^1,\widetilde{a}^1_{j_0}),\le k_{j_0}}(z_0)\}\\
&=\sum_{t=1}^3\min\{n,\nu_{(f^t,\widetilde{a}^t_{j_0}),\le k_{j_0}}(z_0)\}-(2n+1)\min\{1,\nu_{(f^1,\widetilde{a}^1_{j_0}),\le k_{j_0}}(z_0)\}\\
&\ge 2\sum_{j=1}^q\min\{1,\nu_{(f^1,\widetilde{a}^1_j),\le k_j}(z_0)\}+\sum_{t=1}^3\min\{n,\nu_{(f^t,\widetilde{a}^t_{j_0}),\le k_{j_0}}(z_0)\}\\
&\ -(2n+3)\min\{1,\nu_{(f^1,\widetilde{a}^1_{j_0}),\le k_{j_0}}(z_0)\}-2\sum_{t=1}^3\min\{1,\nu_{(f^t,\widetilde{a}^t_{j_0}),>k_{j_0}}(z_0)\}
\end{align*}
 {\bf Case 3.} $z_0 \in \{z \in \C^m: \nu^0_{(f^t,\widetilde{a}^t_{j_0}),> k_{j_0}}>0\}$.
 \begin{align*}
 \nu^0_{\Phi^\alpha}(z_0)& \ge 0 \ge  2\sum_{j=1}^q\min\{1,\nu_{(f^1,\widetilde{a}^1_j),\le k_j}(z_0)\}+\sum_{v=1}^3\min\{n,\nu_{(f^v,\widetilde{a}^v_{j_0}),\le k_{j_0}}(z_0)\}\\
&\ -(2n+3)\min\{1,\nu_{(f^1,\widetilde{a}^1_{j_0}),\le k_{j_0}}(z_0)\}-2\sum_{t=1}^3\min\{1,\nu_{(f^t,\widetilde{a}^t_{j_0}),>k_{j_0}}(z_0)\}
  \end{align*}
Then, from the above three cases it follows that
\begin{align*}
\nu^0_{\Phi^\alpha}(z)&\ge 2\sum_{j=1}^q\min\{1,\nu_{(f^1,\widetilde{a}^1_j),\le k_j}(z_0)\}+\sum_{t=1}^3\min\{n,\nu_{(f^t,\widetilde{a}^t_{j_0}),\le k_{j_0}}(z_0)\}\\
 &-(2n+3)\min\{1,\nu_{(f^1,\widetilde{a}^1_{j_0}),\le k_{j_0}}(z_0)\}-2\sum_{t=1}^3\min\{1,\nu_{(f^t,\widetilde{a}^t_{j_0}),>k_{j_0}}(z_0)\}\\
\end{align*}
for every z outside the analytic set of codimension 2. Integrating both sides of this inequality, we get
\begin{align*}
||N_{\Phi^\alpha}(r)&\ge 2\sum_{j=1}^qN^{[1]}_{(f^1,\widetilde{a}^1_j),\le k_j}(z_0)+\sum_{t=1}^3N^{[n]}_{(f^t,\widetilde{a}^t_{j_0}),\le k_{j_0}}(z_0)\\
 &-(2n+3)N^{[1]}_{(f^1,\widetilde{a}^1_{j_0}),\le k_{j_0}}(z_0)-2\sum_{t=1}^3N^{[1]}_{(f^t,\widetilde{a}^t_{j_0}),>k_{j_0}}(z_0)+o(T(r))\\
\end{align*}
for each $1\le t\le 3$. Hence, the first inequality of lemma is proved.

(b) We now prove the second inequality. By the definition of the Nevanlinna characteristic function, we have
$$
N_{\Phi^\alpha}(r)\le T(r,\Phi^\alpha)+O(1)=N_{\frac{1}{\Phi^\alpha}}(r)+m(r,\Phi^\alpha)+O(1)
$$
We see that a pole of $\Phi^\alpha$ must be zero or pole of $F^{j_0t}_c\ (1\le t\le 3)$. Let $z_0 \not \in \mathcal{V}\cup\mathcal{S}$. There are three cases:

{\bf Case 1.} If $z_0 \in \{z \in \C^m: \nu_{(f^t,\widetilde{a}^t_{j_0}),\le k_{j_0}}^0(z_0)>0\}$, then by (\ref{new1}) we easily see that
$$\nu^\infty_{\Phi^\alpha}(z_0)\le\max_{1\le t\le 3}\nu^\infty_{(f^t,\widetilde{a}^t_{j_0})}+1.$$

{\bf Case 2.} If $z_0 \in \{z \in \C^m: \nu_{(f^t,\widetilde{a}^t_{j_0}),> k_{j_0}}^0(z_0)>0\}$, we rewrite the function $\Phi^\alpha$ as follows
\begin{align*}
\Phi^\alpha(F_{c}^{j_01},F_{c}^{j_02},F_{c}^{j_03}):&=F_{c}^{j_01}F_{c}^{j_02}F_{c}^{j_03}\left | \begin {array}{cccc}
1&1&1\\
F_{j_0}^{c1}&F_{j_0}^{c2} &F_{j_0}^{c3})\\
\mathcal {D}^{\alpha}(F_{j_0}^{c1}) &\mathcal {D}^{\alpha}(F_{j_0}^{c2}) &\mathcal {D}^{\alpha}(F_{j_0}^{c3})\\
\end {array}
\right|\\
&=F^{j_01}_c(F^{j_02}_c-F^{j_03}_c)D^\alpha(F^{c1}_{j_0})+F^{j_02}_c(F^{j_03}_c-F^{j_01}_c)D^\alpha(F^{c2}_{j_0})\\
&+F^{j_03}_c(F^{j_01}_c-F^{j_02}_c)D^\alpha(F^{c3}_{j_0})
\end{align*}
It is easy to see that $\nu_{\Phi^\alpha}^\infty(z_0)\le \max_{1\le t\le3}\{\nu^\infty_{F^{j_0t}_cD^\alpha(F^{ct}_{j_0})}(z_0)\}\le |\alpha|=1$

{\bf Case 3.} If $z_0 \in \{z \in \C^m: \nu_{(f^t,a^t_c)}>0\}$ then $\nu_{\Phi^\alpha}^\infty(z_0) \le \sum_{t=1}^3\nu_{F^{j_0t}_c}^\infty(z_0)$.

Thus, every $z \not \in\mathcal{V}\cup\mathcal{S}$, we have $\nu_{\Phi^\alpha}^\infty(z_0) \le \sum_{t=1}^3\nu_{F^{j_0t}_c}^\infty(z_0)+\sum_{t=1}^3\min\{\nu_{(f^t,\widetilde{a}^t_{j_0}),>k_{j_0}},1\}$. Therefore, we have
\begin{align*}
||N_{\frac{1}{\Phi^\alpha}}(r)&\le \sum_{t=1}^{3}N_{\frac{1}{F_c^{j_0t}}}(r)+\sum_{t=1}^3N^{[1]}_{(f^t,\widetilde{a}^t_{j_0}),>k_{j_0}}(r)+o(T(r))
\end{align*}
By the logarithmic derivative lemma (Lemma $\ref{2.1}$), we have
\begin{align*}
|| m(r,\Phi^\alpha)&\le \sum_{t=1}^{3}m(r,F_{c}^{j_0t})+O\biggl(\sum m\biggl(r,\dfrac{{\mathcal D}^{\alpha^i}(F_{j_0}^{ct})}{F_{j_0}^{ct}}\biggl)\biggl)+O(1)\\
&\le \sum_{t=1}^{3}m(r,F_{c}^{j_0t})+\sum_{t=1}^{3}o(T(r,F_{j_0}^{ct}))+O(1)\\
&=\sum_{t=1}^{3}m(r,F_{c}^{j_0t})+o(T(r)).
\end{align*}
This implies that
\begin{align*}
\biggl|\biggl|N_{\Phi^\alpha}(r)&\le \sum_{t=1}^{3}T(r,F_c^{j_0t})+o(T(r))\le \sum_{t=1}^{3}T_{f^t}(r)+\sum_{t=1}^3N^{[1]}_{(f^t,\widetilde{a}^t_{j_0}),>k_{j_0}}(r)+o(T(r))\\
&\le T(r)+\sum_{t=1}^3N^{[1]}_{(f^t,\widetilde{a}^t_{j_0}),>k_{j_0}}(r)+o(T(r))\quad \hfill\square
\end{align*}

\vskip0.2cm
\noindent
{\bf 3.2. Proof of Theorem \ref{1.2}.}\
 Assume that
$$\sum_{i=1}^q\dfrac{1}{k_i+1} < \dfrac{q-3n+2}{2q-5n+10}\biggl(\dfrac{2(2q+n-3)}{3n(n+1)}-3\biggl).$$

Denote by $\mathcal{Q}$ be the set of all indices $j\in \{1,..,q\}$ satisfying the following: there exist $c\in \mathcal{C}$ and   $\alpha=(\alpha_1,\ldots ,\alpha_m)\in \Z_+^m $ with $|\alpha|= 1$ such that $\Phi^\alpha (F_{c}^{j1},F_{c}^{j2},F_{c}^{j3}) \not \equiv 0$. We put $p=\sharp \mathcal{Q}$.

Suppose that $p \ge q-3n+2$. Without loss  of generality, we may assume that $1,...,q-3n+2 \in \mathcal{Q}$. Then by Lemma \ref{4.7}, for $j \in \mathcal{Q}$, $1\le t\le 3$, we have
%Suppose that there exist two indices $i,j\in\{1,\ldots ,q\}$ and  $\alpha=(\alpha_1,\alpha_2)\in (\Z_+^n)^2 $ with $|\alpha|\le 1$ such that $\Phi^\alpha (F_{j}^{i1},F_{j}^{i2},F_{j}^{i3}) \not \equiv 0$. By Lemma \ref{4.7}, we have
\begin{align*}
 T(r) &\ge2\sum_{j=1}^qN^{[1]}_{(f^v,\widetilde{a}^v_j),\le k_j}(r)+\sum_{t=1}^3N^{[n]}_{(f^t,\widetilde{a}^t_{i}),\le k_{i}}(r)\\
 &-(2n+3)N^{[1]}_{(f^v,\widetilde{a}^v_{i}),\le k_{i}}(r)-3\sum_{t=1}^3N^{[1]}_{(f^t,\widetilde{a}^t_{i}),>k_{i}}(r)+o(T(r)).
\end{align*}
By summing up both side of above inequality over $1\le j\le q-3n+2$ and $1\le t \le 3$, we have
\begin{align*}
\parallel 3(q-3n+2)T(r)&\ge 2(q-3n+2)\sum_{t=1}^3 \sum_{i=1}^qN_{(f^t,\widetilde{a}_i^t),\leq k_i}^{[1]}(r)+3\sum_{t=1}^{3}\sum_{i=1}^{q-3n+2}N^{[n]}_{(f^t,\widetilde{a}^t_{i}),\le k_{i}}(r)\\
&-(2n+3)\sum_{t=1}^3\sum_{i =1}^{q-3n+2}N^{[1]}_{(f^t,\widetilde{a}^t_{i}),\le k_{i}}(r)-9\sum_{t=1}^3\sum_{i=1}^{q-3n+2}N^{[1]}_{(f^t,\widetilde{a}^t_{i}),>k_{i}}(r)\\
&= \sum_{t=1}^3\biggl((2q-8n+1) \sum_{i=1}^{q-3n+2}N_{(f^t,\widetilde{a}_i^t),\leq k_i}^{[1]}(r)+3\sum_{i=1}^{q-3n+2}N^{[n]}_{(f^t,\widetilde{a}^t_{i}),\le k_{i}}(r)\\
&+(2q-6n+4)\sum_{i =q-3n+3}^qN^{[1]}_{(f^t,\widetilde{a}^t_{i}),\le k_{i}}(r)-9\sum_{i=1}^{q-3n+2}N^{[1]}_{(f^t,\widetilde{a}^t_{i}),>k_{i}}(r)\biggl)\\
&= \sum_{t=1}^3\biggl((2q-8n+1) \sum_{i=1}^{q-3n+2}N_{(f^t,\widetilde{a}_i^t)}^{[1]}(r)+3\sum_{i=1}^{q-3n+2}(N^{[n]}_{(f^t,a^t_{i})}(r)-N^{[n]}_{(f^t,\widetilde{a}^t_{i}),>k_i}(r))\\
&+(2q-6n+4)\sum_{i =q-3n+3}^q(N^{[1]}_{(f^t,\widetilde{a}^t_{i})}(r)-N^{[1]}_{(f^t,\widetilde{a}^t_{i}),>k_i}(r))\\
&\ -(2q-8n+10)\sum_{i=1}^{q-3n+2}N^{[1]}_{(f^t,\widetilde{a}^t_{i}),>k_{i}}(r)\biggl)\\
&\ge \sum_{t=1}^{3}\biggl(\dfrac{2q-5n+1}{n} \sum_{i=1}^{q-3n+2}N_{(f^t,\widetilde{a}_i^t)}^{[n]}(r)+\dfrac{2q-6n+4}{n}\sum_{i =q-3n+3}^qN^{[n]}_{(f^t,\widetilde{a}^t_{i})}(r)\\
&-\sum_{i=1}^{q-3n+2}\dfrac{2q-5n+10}{k_i+1}N_{(f^t,\widetilde{a}^t_{i})}(r)-\sum_{i=q-3n+3}^q\dfrac{2q-6n+4}{k_i+1}N_{(f^t,\widetilde{a}^t_i)}(r)\biggl)\\
&\ge \sum_{t=1}^{3}\biggl(\dfrac{2q-5n+1}{n} \sum_{i=1}^{q-3n+2}N_{(f^t,a_i^t)}^{[n]}(r)+\dfrac{2q-6n+4}{n}\sum_{i =q-3n+3}^qN^{[n]}_{(f^t,a^t_{i})}(r)
\end{align*}
\begin{align*}
&-\sum_{i=1}^{q-3n+2}\dfrac{2q-5n+10}{k_i+1}N_{(f^t,a^t_{i})}(r)-\sum_{i=q-3n+3}^q\dfrac{2q-6n+4}{k_i+1}N_{(f^t,a^t_i)}(r)\biggl)+o(T(r))\\
\end{align*}
%Note that, here we use the assumption that $ q\ge 3n+3>4n-\dfrac{1}{2},$ then $2q-8n+1 \ge 0$.
%, and hence $(2q-8n+1)N^{[1]}_{(f^t,a^t_i),\le k_i}\ge (\dfrac{2q-8n+1}{n}+3)N^{[n]}_{(f^t,a^t_j),\le k_j})$

On the other hand, by theorem \ref{2.3}, we have
\begin{align*}
|| 3(q-3n+2)T(r)&\ge \biggl(\dfrac{2q-5n+1}{n}.\dfrac{2(q-3n+2)}{3(n+1)} +\dfrac{2q-6n+4}{n}.\dfrac{2(3n-2)}{3(n+1)}\\
&\quad -\sum_{i=1}^{q-3n+2}\dfrac{2q-5n+10}{k_i+1}-\sum_{i=q-3n+3}^{q}\dfrac{2q-6n+4}{k_i+1}\biggl)T(r)+o(T(r))\\
&\ge \biggl(\dfrac{(2q-6n+4)(2q+n-3)}{3n(n+1)}-(2q-5n+10)\sum_{i=1}^{q}\dfrac{1}{k_i+1}\biggl)T(r)+o(T(r))\\
\end{align*}
Letting $r\longrightarrow +\infty$, we get
\begin{align*}
&3(q-3n+2)\ge\dfrac{(2q-6n+4)(2q+n-3)}{3n(n+1)}-(2q-5n+10)\sum_{i=1}^{q}\dfrac{1}{k_i+1}\\
&\text{ i.e., } \sum_{i=1}^q\dfrac{1}{k_i+1} \ge \dfrac{q-3n+2}{2q-5n+10}\biggl(\dfrac{2(2q+n-3)}{3n(n+1)}-3\biggl)
\end{align*}
 This is a contradiction.

Then $\sharp\mathcal{Q}\le q-3n+1$. Without loss of generality, we may assume that $1,2,...,3n-1 \not \in \mathcal{Q}$.
 This mean that
$$\Phi^\alpha (F_{c}^{j1},F_{c}^{j2},F_{c}^{j3}) \equiv 0,$$
for all $c \in \mathcal C$, $\alpha=(\alpha_1,\ldots ,\alpha_m)\ \text { with }|\alpha|= 1$. By the density of $\mathcal C$ in $\C^q$, the above equality holds for all $c \in \C^q$, and $|\alpha|= 1$.  For each $i \in \{1,...,3n-1\}$, chosing $c_i=(0,...,0,\underset{i-th}{1},0,...,0)$ we have
$$\Phi^\alpha (F_{c_i}^{j1},F_{c_i}^{j2},F_{c_i}^{j3}) \equiv 0\ \forall |\alpha|=1.$$
Then by Lemma \ref{4.4}, there exists a constant $\lambda$ such that
$$F_{c_i}^{j1}=\lambda F_{c_i}^{j2},F_{c_i}^{j2}=\lambda F_{c_i}^{j3}, \text { or } F_{c_i}^{j3}=\lambda F_{c_i}^{j1}.$$
For instance, we assume that $F_{c_i}^{j1}=\lambda F_{c_i}^{j2}$. We will show that $ \lambda=1.$

Indeed, suppose that $\lambda \ne 1$, we have
\begin{align*}
0=T(r,\dfrac{F^{j1}_{c_i}}{F^{j2}_{c_i}})&\ge N^{[1]}_{F^{j1}_{c_i}-F^{j2}_{c_i}}(r)\ge \sum_{\underset{v\ne i,j}{v=1}}^qN^{[1]}_{(f^1,\widetilde{a}^1_v),\le k_v}-\sum_{t=1}^2\sum_{v=i,j}N^{[1]}_{(f^t,\widetilde{a}^t_v),> k_v}(r)\\
&\ge \dfrac{1}{2}\sum_{\underset{v\ne i,j}{v=1}}^q\biggl(N^{[1]}_{(f^1,\widetilde{a}^1_v),\le k_v}+N^{[1]}_{(f^2,\widetilde{a}^2_v),\le k_v}\biggl)-\sum_{t=1}^2\sum_{v=i,j}N^{[1]}_{(f^t,\widetilde{a}^t_v),> k_v}(r)\\
&\ge \dfrac{1}{2}\sum_{\underset{v\ne i,j}{v=1}}^q\biggl(N^{[1]}_{(f^1,\widetilde{a}^1_v)}+N^{[1]}_{(f^2,\widetilde{a}^2_v)}-N^{[1]}_{(f^1,\widetilde{a}^1_v),> k_v}-N^{[1]}_{(f^2,\widetilde{a}^2_v),>k_v}\biggl)\\
&\ -\sum_{t=1}^2\sum_{v=i,j}N^{[1]}_{(f^t,\widetilde{a}^t_v),> k_v}(r)\\
&\ge \dfrac{1}{2}\sum_{v\ne i}\biggl(\dfrac{1}{n}(N^{[n]}_{(f^1,a^1_v)}+N^{[n]}_{(f^2,a^2_v)})-\dfrac{1}{k_v+1}(N_{(f^1,a^1_v) }+N_{(f^2,a^2_v)})\biggl)\\
&\ -\sum_{t=1}^2\sum_{v=i,j}\dfrac{1}{k_v+1}N_{(f^t,a^t_v)}(r)+o(T(r))\\
&\ge\biggl( \dfrac{q-2}{3n(n+1)}-\dfrac{1}{2}\sum_{\underset{v\ne i,j}{v=1}}^q\dfrac{1}{k_v+1}-\sum_{v=i,j}\dfrac{1}{k_v+1}\biggl)T(r)+o(T(r))\\
&\ge\biggl( \dfrac{q-2}{3n(n+1)}-\sum_{v=1}^q\dfrac{1}{k_v+1}\biggl)T(r)+o(T(r))\\
\end{align*}
Thus, $\sum_{v=1}^q\dfrac{1}{k_v+1}\ge \dfrac{q-2}{3n(n+1)}$.
%Since $F_{i}^{j1}=F_{i}^{j2}$ on the set $\bigcup_{v\ne i}\{z : (f^t,a^t_v)(z)=0\},$ we have that
%$F_{i}^{j1}=F_{i}^{j2}=0$ on the set $\bigcup_{v\ne i}\{z : (f^t,a^t_v)(z)=0\}.$ Hence
%$\bigcup_{v\ne i}\{z : (f^t,a^t_v)(z)=0\} \subset \{z : (f^t,a^t_i)(z)=0\}.$ It follows that $\{z : (f^t,a^t_v)(z)=0\}=\emptyset \ (v\ne i,j).$
%We obtain that
%$$\parallel\dfrac{2(q-2)}{3(n+1)}T_{f^t}(r) \le \sum_{v\ne i,v\ne j}N_{(f^t,a^t_v),\le k_v}^{[n]}(r)+o(T_{f^t}(r))=o(T_{f^t}(r)).$$
This is a contradiction. Thus $ \lambda =1\ (1\le i<j\le q).$

Define
$$I_1= \{j\in \{2,\ldots ,3n-1\}: F_{1}^{j1}=F_{1}^{j2}\},$$
$$I_2= \{j\in \{2,\ldots ,3n-1\}: F_{1}^{j2}=F_{1}^{j3}\},$$
$$I_3= \{j\in \{2,\ldots ,3n-1\}: F_{1}^{j3}=F_{1}^{j1}\}.$$
Since $\sharp (I_1\cup I_2\cup I_3)=\sharp \{2,\ldots ,3n-1\}=3n-2$, there exists $1\le v\le 3$ such that $\sharp \ I_v \ge n$.
Without loss of generality, we may assume that $\sharp \ I_1 \ge n$. This implies that
$$\dfrac{(f^1,\widetilde{a}^1_1)}{(f^2,\widetilde{a}^2_1)}=\dfrac{(f^1,\widetilde{a}^1_j)}{(f^2,\widetilde{a}^2_j)}\ \forall j\in I_1.$$
The theorem is proved.\hfill$\square$

\end{document}